\documentclass[12pt,oneside,english]{amsart}
\usepackage[latin1]{inputenc}
\usepackage{amssymb}

\makeatletter
\newtheorem{theorem}{Theorem}
\newtheorem{lemma}{Lemma}

\newtheorem{conjecture}{Conjecture}
\newtheorem{problem}{Problem}

\numberwithin{equation}{section}
\usepackage{babel}

\makeatother
\begin{document}
\baselineskip=17pt

\title[Infinite sets of generators of primes ]{Infinite sets of generators of primes based on the Rowland idea and  conjectures concerning twin primes}

\author{Vladimir Shevelev}
\address{Departments of Mathematics \\Ben-Gurion University of the
 Negev\\Beer-Sheva 84105, Israel. e-mail:shevelev@bgu.ac.il}


\begin{abstract}
Using the Rowland idea, we find two infinite sets of generators of primes. We also pose some conjectures concerning
twin primes.
\end{abstract}

\maketitle

\section{Introduction }

In 2008, Rowland [1], using very elementary tools, discovered a very interesting fact.
\begin{theorem}\label{1}
Let $a(1)=7$ and, for $n\geq2,$
$$a(n)=a(n-1)+\gcd (n, \enskip a(n-1)). $$
Then, for every $n\geq2,$ the difference $a(n)-a(n-1)$ is 1 or prime.
 \end{theorem}
   Theorem 1 resembles a celebrated chess problem of Richard Reti with one white and one black pawns. Earlier (2006)
   A. Adamchuk obtained a more complicated one-prime generators of primes of special kinds (see sequences A118679 and A120293 in [3]) and some later (2008) B. Cloitre proposed a conjugate to the Rowland's one-prime generators  with lcm instead of gcd (see sequences A135506, A135508 in [3], and the details of his approach one can see in a very interesting review paper [2]). If to ignore a possible practical applications, I think that the Rowland one-prime sequence is the most attractive in view of it has very many 1's. A meeting of a large prime after, say, more than million 1's, produces
   a very strong impression, maybe, similar to a small "big bang". The average  value $\overline{a}(N)$ of $a(n)-a(n-1)$ on [2,N] is $(a(N)-7)/N.$ Therefore, from the Rowland results  it follows that $\lim\sup \overline{a}(N)=3.$   If to compare the ubiquitous 1's of the Rowland sequence with the atoms of hydrogen in the universe, then it is interesting to note that, according to Einstein's theory, an "open" universe has a infinite volume and an infinite number of hydrogen atoms and the universe is at the critical density which is found to be 3 hydrogen atoms per cubic meter. This mystical coincidence, in my opinion, is very interesting! \newline
\section{Infinite sets of generators of primes based on the Rowland idea}
 The following theorem gives an infinite set of Rowland-like generators of primes (cf Theorem 1)
 \begin{theorem}\label{2}
Let, for  $m\geq2,\enskip c(m-1)=t$ such that $\gcd(m,\enskip t)$ is a prime. Denote $\gcd(m,\enskip t)$ by $p.$ Let  \enskip$ t+p=3m.$ If, for $n\geq m,$
$$c(n)=c(n-1)+\gcd (n, \enskip c(n-1)), $$
 then, for every $n\geq m,$ the difference $c(n)-c(n-1)$ is 1 or prime.
 \end{theorem}
 \indent Note that we easily could choose $m$ and $t$ such that the first difference is an arbitrary prime $P.$
 E.g., if $m=P$ and $t=2P,$ then $c(m)=2P+\gcd(2P,P)=3P; \enskip c(m)-c(m-1)=P $ (and conditions of Theorem 2 are, evidently, satisfied; cf. A167170 in [3]).\newline
 Denoting $\{c_P(n)\}$ the sequence of Theorem 2 with initial condition $c_P(P-1)=2P,$ note that, for $n>=3,$ sequences  $\{c_2(n)\}$ and $\{a(n)\}$ from Theorem 1 coincide. Let us construct the sequence of primes for which the generators of primes $\{c_P(n)\}$ are different. Put $P_1=2.$ Note that, the least $P>2$ for which $c_2(P-1)$ does not equal to $2P$ is 7; therefore put $P_2=7.$ The least $P>7,$ for which $c_2(P-1)$ and $c_7(P-1)$ are different from $2P,$ is 17, therefore, put $P_3=17,$ etc. We obtain sequence (cf. A167168 [3]):
 $$ 2, 7, 17, 19, 37, 43, 53, ...$$
\indent Proof of Theorem 2 is based on the following version of Lemma 1[1].
\begin{lemma}
Let
 $$c(n)=c(n-1)+\gcd (n,\enskip c(n-1)).$$
Put
$$h(n)= c(n)-c(n-1).$$
Let, for $n_1\geq2, \enskip c(n_1) =3n_1.$ \newline
Finally, let $n_2$ be the smallest integer greater than $n_1$ such that $h(n_2)>1.$\newline
Then $c(n_2)=3n_2$ and $h(n_2)$ is prime.
\end{lemma}

Quite analogously we obtain the following result
\begin{theorem}\label{4}
Let, for $m\geq4,\enskip c(m-1)=t$ such that $\gcd(m,\enskip t)$ is a prime. Denote $\gcd(m,\enskip t)$ by $p.$ Let  \enskip$ t+p=2m.$ If, for $n\geq m,$
$$c(n)=c(n-1)+\gcd (n, \enskip c(n-1)), $$
 then, for every $n\geq m,$ the difference $c(n)-c(n-1)$ is 1 or prime.
 \end{theorem}\newpage
 \section{Conjectures concerning twin primes }
 \begin{conjecture}\label{1}
Let $c(1)=2$ and for $n\geq 2,$
$$c(n)=c(n-1)+\begin{cases}\gcd (n, \enskip c(n-1)),\enskip if \;\; n\enskip is\enskip even
\\\gcd (n-2,\enskip c(n-1)),\enskip if\;\; n\enskip is\enskip odd\end{cases}. $$
 Then every record (more than 3) of the value of difference $c(n)-c(n-1)$ is greater of twin primes.
 \end{conjecture}
 The first records are (cf. sequence A166945 in [3])
 \begin{equation}\label{3.1}
  7,13,43,139,313,661,1321,2659,5459,10891,22039,...
  \end{equation}
For "conjugate" sequence, which is defined by the recurrence
$$l(1)=2, \enskip l(n)=l(n-1)+\begin{cases}\gcd (n,\enskip l(n-1)), \enskip if \;\; n\enskip is\enskip odd,
\\\gcd (n-2,\enskip l(n-1)),\enskip if\;\; n\enskip is\enskip even.\end{cases} $$
we pose a stronger conjecture.
 \begin{conjecture}\label{2}
1) For $n\geq2,$ every difference $l(n)-l(n-1)$ is 1 or prime;\newline
 2) Every record of the value of difference $l(n)-l(n-1)$ is greater of twin primes.
 \end{conjecture}
 Here the first records are (cf. sequence A167495 in [3])
 \begin{equation}\label{3.2}
  5,13,31,61,139,283,571,1353,2911,4651,9343,19141,...
   \end{equation}
   Conjectures 1 and 2 are connected with the following general problem.
 \begin{problem}
 To find a generator of the lower (higher) primes of twin pairs or products of lower and higher primes of twin pairs.
 \end{problem}
 A solution of this problem in a Rowland manner (cf. our sequence A167054 in [3]) could, maybe, throw light on the infinity of twin primes.\newline
 \section{Conjectures concerning twin primes using function $s(m,n)$}
 Denote $\mathrm{P}$ the set of primes, $\mathrm{P}_1=\mathrm{P}\cup \{1\}.$
For positive integers $m,n,$ denote $s(m,n)$ 
 the following function $\mathbb{N}\times\mathbb{N}\mapsto\mathrm{P}_1:$
 \begin{equation}\label{4.1}
  s(m,n)=\begin{cases}1,\enskip if \;\; \gcd(m,n)=1,
\\\max{(p\in \mathrm{P}, \enskip p|\gcd(m,n))},\;\; otherwise.\end{cases}
   \end{equation}
Of course, this function is senseless with the point of view a constructing of any generator of primes, but is very
  useful for the formulation of conjectures like Conjectures 1 and 2. We illustrate this on the following conjectures
  close to  Conjectures 1, 2.\newpage
  \begin{conjecture}\label{3}
Let $c(1)=2$ and for $n\geq 2,$
$$c(n)=c(n-1)+\begin{cases} s(n, \enskip c(n-1)),\enskip if \;\; n\enskip is\enskip even
\\s(n-2,\enskip c(n-1)),\enskip if\;\; n\enskip is\enskip odd.\end{cases} $$
 Then every record (more than 3) of the value of difference $c(n)-c(n-1)$ is greater of twin primes.
 \end{conjecture}
 Instead of (3.1), here we have sequence
 \begin{equation}\label{4.2}
  7,13,31,61,151,313,661,1321,2659,5459,10891,22039,...
  \end{equation}
  (and it is natural to conjecture that, beginning with 313 sequences (3.1) and (4.2) coincide). The evident advantage
  of such formulation of Conjecture 3 for the research and the verification consists of the fact that we should not verify the terms of sequence (4.2) on the primality: all of them are primes by the definition of the sequence.
  Quite analogously we formulate the "conjugate" conjecture.
   \begin{conjecture}\label{4} Let $l(1)=2$ and for $n\geq 2,$
   $$ l(n)=l(n-1)+\begin{cases}s(n,\enskip l(n-1)), \enskip if \;\; n\enskip is\enskip odd,
\\s(n-2,\enskip l(n-1)),\enskip if\;\; n\enskip is\enskip even.\end{cases} $$
Then every record (more than 3) of the value of difference $l(n)-l(n-1)$ is greater of twin primes.
 \end{conjecture}
\;\;\;\;\;\;\;\;

\;\;\;\;\;\;\;\;

\end{document}